\documentclass[a4paper]{article}

\usepackage[latin1]{inputenc} %
\usepackage[T1]{fontenc} %
\usepackage{RR}
\usepackage{hyperref}

\usepackage{graphicx,amsmath,amssymb}
\usepackage{subfigure}
\usepackage{epsfig}
\usepackage{verbatim}
\usepackage{algorithm}
\usepackage{algorithmic}
\usepackage{color}
\usepackage{multirow}
\usepackage{latexsym}

\RRdate{September 2010}

\RRauthor{%
F. Grognard\thanks{INRIA Sophia Antipolis-M\'editerran\'ee, France, Frederic.Grognard@inria.fr}
\and
A.R. Akhmetzhanov\thanks{INRIA Sophia Antipolis-M\'editerran\'ee, France, akhmetzhanov@gmail.com}
\and
P. Masci\thanks{INRIA Sophia Antipolis-M\'editerran\'ee, France, Pierre.Masci@inria.fr}
\and
O. Bernard\thanks{INRIA Sophia Antipolis-M\'editerran\'ee, France, Olivier.Bernard@inria.fr}
}

\authorhead{Grognard \& Akhmetzhanov \& Masci \& Bernard}

\RRtitle{Optimisation de la production de biomasse par un photobioréacteur en utilisant la lumière naturelle}
\RRetitle{Optimization of a photobioreactor biomass production using natural light}
\titlehead{Optimization of a photobioreactor biomass production using natural light}

\RRresume{Dans ce rapport, nous abordons la question de l'optimisation de la productivité à long terme de la biomasse microalgale dans le cadre d'une
production en photobioréacteur sous l'influence du cyle jour/nuit. Pour cela, nous proposons un modèle simple de bioréacteur
représentant l'atténuation de la lumière dans le réacteur due à l'auto-ombrage de la biomasse. Nous obtenons une loi de commande qui utilise le taux de
dilution comme contrôle et optimise la productivité sur une seule journée par
l'application du principe du maximum de Pontryagin. Une contrainte importante à la solution obtenue est que la biomasse dans le réacteur devrait être au même niveau au début et à la fin de la journée pour que le même contrôle puisse être appliqué tous les jours afin d'optimiser la productivité à long terme. Plusieurs scénarios sont possibles en fonction des paramètres du modèle de croissance de la micro-algues et de la valeur maximale admissible du taux de dilution:
commande bang-bang ou bang-singulière-bang  ou, si le taux de croissance des algues est très fort en présence de la lumière, dilution maximale constante. Un diagramme de bifurcation est présenté pour illustrer pour quelles valeurs des paramètres ces différents comportements se produisent.
}

\RRabstract{
We address the question of optimization of the biomass long term productivity in the framework of microalgal biomass
production in photobioreactors under the influence of day/night cycles. For that, we propose a simple bioreactor model
accounting for light attenuation in the reactor due to biomass density and obtain the control law that optimizes productivity over a single day through
the application of Pontryagin's maximum principle, with the dilution rate being the control. An important constraint on the obtained solution is that the biomass in the reactor should be at the same level at the beginning and at the end of the day so that the same control can be applied everyday and optimizes the long term productivity. Several scenarios are possible depending on the microalgae's strain parameters and the maximal admissible value of the dilution rate:
bang-bang or bang-singular-bang control or, if the growth rate of the algae is very strong in the presence of light, constant maximal dilution. A bifurcation diagram is presented to illustrate for which values of the parameters these different behaviors occur.
}

\RRmotcle{Contrôle optimal, systèmes biologiques, modélisation}
\RRkeyword{Optimal control; Biological systems; Modeling}

\RRprojet{Comore}  
\RRtheme{\THBio} 
\URSophia 

\begin{document}

\RRNo{7378}

\makeRR

\section{Introduction}

Microalgae  have recently received more and more attention in the frameworks of CO$_2$ fixation and renewable energy \cite{Huntley2007, Chisti2007}. Their high actual
photosynthetic yield compared to terrestrial plants (whose growth is limited by CO$_2$
availability) leads to large potential algal biomass productions in photobioreactors of several tens of tons
per hectare and per year \cite{Chisti2007}.

The objective of this paper is to develop an optimal control law that would
maximize the photobioreactor yield, while taking into account that the light
source ({\it i.e} the primary energy source) that will be used is the natural light.
The light source is therefore
periodic with a light phase (day) and a dark phase (night). In addition to
this time-varying periodic light source, we will take the auto-shading in the
photobioreactor into account: the pigment concentration (mainly chlorophyll)
affects the light distribution and thus the biological activity within the reactor. As a consequence, for a too
high biomass, light in the photobioreactor is strongly attenuated and growth
is low.

It is therefore necessary to develop a model that takes both  features
into account in order to develop the control law, where the substrate
concentration in the input (marginally) and the dilution rate (mainly) will be
used. This model should not be too complicated in order to be tractable and should
present the main features of the process. Since we want to develop a control strategy that will be used on the
long run, we could choose an infinite time-horizon measure of the
yield. However, we rather took advantage of the observation that, in the
absence of a discount rate in the cost functional, the control should be
identical everyday and force the state of the system to be identical at the
beginning of the day and 24 hours later. We therefore opted for optimizing a
cost over one day with the constraint that the initial and terminal state
should be identical.

The paper is structured as follows: first, we present the model dealing with
both substrate limitation, light attenuation and light periodicity; then
biomass productivity optimization is presented in a constant light
environment. The solution to the periodic light problem is then
presented. Finally, numerical results are presented
with a bifurcation analysis.

\section{A photobioreactor model with light attenuation}

Micro-algae growth in a photobioreactor is often modelled through one of two
models, the Monod model \cite{Monod1942} or the Droop Model
\cite{Droop1968}. The latter is more accurate as it separates the process of substrate uptake and growth of the microalgae. The former gives a reasonable representation of reality by coupling growth and uptake, and is more convenient for building control laws since it is simpler. For sake of simplicity we will introduce the problem with the Monod model, but the presented results are similar with the Droop model when considering the working modes where nutrients are not limiting growth. The Monod model writes:
\begin{equation}\label{monod}
\left\{
\begin{array}{lll}
  \frac{ds}{d\tau}& =& D(s_{in}-s)-k\nu(s)x\\
  \frac{dx}{d\tau}& =& \nu(s)x-Dx
\end{array}
\right.
\end{equation}
where $s$ and $x$ are the substrate and biomass concentrations in the medium,
while $D$ is the dilution rate, $s_{in}$ is the substrate input concentration and
$k$ is the substrate/biomass yield coefficient. We will depart from this model in two
directions.  First, we introduce respiration
by the microalgae: contrary to photosynthesis, this phenomeneon takes place
with or without light; from a carbon point of view, it converts biomass into
carbon dioxyde, so that we represent it as a $-\rho x$ term in the biomass
dynamics. Secondly, under the hypothesis of an horizontal planar photobioreactor (or raceway) with  vertical incoming light,
we represent light attenuation following an exponential
Beer-Lambert law where the attenuation at some depth $z$ comes from the total
biomass $xz$ per surface unit contained in the layer of depth $[0,z]$:
\begin{equation} \label{eq:Iz}
I(xz) = I_0 e^{-a x z}
\end{equation}
where $I_0$ is the incident light and $a$ is a light attenuation coefficient.
In microalgae, as we proposed in (\ref{eq:Iz})
chlorophyll is mostly the cause of this shadow effect and, in
model (\ref{monod}), it is best represented by a fixed portion of the biomass \cite{Bernard2009}. Finally,
the light source variation will be introduced by
taking a time-varying incident light $I_0(\tau)$. With such an hypothesis on the
light intensity that reaches depth $z$, growth rates vary with depth: in the upper part of the reactor, higher light causes higher growth than in the bottom part. Supposing that light
attenuation directly affects the maximum growth rate \cite{Huismann2002}, the growth rate for a given depth $z$ can then be written as
$$
\begin{array}{l}
\displaystyle \nu_z(s,I(xz,\tau)) = \frac{ \tilde\nu I(xz,\tau)}{I(xz,\tau)+K_I}\frac{s}{s+K_s},\\
\hspace*{4cm}\mbox{with }I(xz,\tau) = I_0(\tau)e^{-axz}
\end{array}
$$
Then, we can compute the mean growth rate in the reactor:
\[
\nu(s,I_0(\tau),x) = \frac{1}{L} \int_0^L \nu_z(s,I(xz,\tau)) dz
\]
where $L$ is the depth of the reactor and where we have supposed that, even
though the growth rate is not homogeneous in the reactor due to the light
attenuation, the concentrations of $s$ and $x$ are kept homogeneous through
continuous reactor stirring. It is this average growth rate that will be used in the lumped
model that we develop.  We then have:
$$
\begin{array}{lll}
\nu(s,I_0(\tau),x)&=& \displaystyle\frac{\tilde\nu}L\int_0^L\frac{I_0(\tau)e^{-a xz}}{I_0(\tau)e^{-a xz}+K_I}\>dz\frac{s}{s+K_s}\\
&=& \displaystyle\frac{\tilde\nu}{a xL}\ln\left(\!\frac{I_0(\tau)+K_I}{I_0(\tau)e^{-a xL}+K_I}\!\right)\!\frac{s}{s+K_s}
\end{array}
$$

The system for which we want to build an optimal controller is therefore
\begin{equation}\label{our_monod}
\left\{
\begin{array}{lll}
  \frac{ds}{d\tau}& =& D(s_{in}-s)-k \frac{\tilde\nu}{a xL}\ln\left(\!\frac{I_0(\tau)+K_I}{I_0(\tau)e^{-a xL}+K_I}\!\right)\!\frac{s}{s+K_s} x\\
  \frac{dx}{d\tau}& =& \frac{\tilde\nu}{a xL}\ln\left(\!\frac{I_0(\tau)+K_I}{I_0(\tau)e^{-a xL}+K_I}\!\right)\!\frac{s}{s+K_s}x-\rho x-Dx
\end{array}
\right.
\end{equation}
However, since we want to maximize the productivity, it seems clear that the
larger $s$ the better, large values of $s$ translating into large growth rates. The control $s_{in}$ should then always be kept very
large so as to always keep the substrate in the region where
$\frac{s}{s+K_s}\approx 1$. We can then concentrate on the reduced model
\begin{equation}\label{reduced}
  \frac{dx}{d\tau} = \frac{\tilde\nu}{a xL}\ln\frac{I_0(\tau)+K_I}{I_0(\tau)e^{-a xL}+K_I}x-\rho x-Dx
\end{equation}
which then encompasses all the relevant dynamics for the control problem.

In order to more precisely determine the model, we should now indicate what
the varying light will be like. Classically, it is considered that daylight
varies as the square of a sinusoidal function so that
\[
I_0(\tau)=\left(\max\left\{\sin\left(\frac{2\pi \tau}{{\mathcal T}}\right),0\right\}\right)^2
\]
where ${\mathcal T}$ is the length of the day. The introduction of such a varying light
would however render the computations analytically untractable. Therefore, we approximate the light source by a step
function:
$$
I_0(\tau)=\left\{\begin{array}{lllll}
  \bar I_0,&\  & 0\le \tau < \bar {\mathcal T}&\mbox{\quad ---\quad light  phase}\\
  0,&\  & \bar {\mathcal T}\le \tau< {\mathcal T}&\mbox{\quad ---\quad dark phase}
\end{array}\right.
$$
In a model where the time-unit is the day, ${\mathcal T}$ will be equal to $1$. In the
following, we will consider $\bar {\mathcal T}=\frac{{\mathcal T}}{2}$, but this quantity obviously
depends on the time of the year.

Finally, we consider a last simplification to the model: instead of
considering that the biomass growth in the presence of light has the form
$\frac{\tilde\nu}{a L}\ln\frac{\bar I_0+K_I}{\bar I_0 e^{-a xL}+K_I}$, which is
an increasing and bounded function, we replace it with another increasing
bounded function $\frac{\bar \nu x}{k+x}$ and
obtain the model
\[
  \frac{ dx}{d\tau} = \frac{\nu(\tau) x}{\kappa+x}-\rho x-Dx
\]
where $\nu(\tau)=\bar \nu$ during the light phase and $0$ at night. It is possible to show that this simplified model is a good  numerical approximation of the original model.

\section{Productivity optimization}

The productivity problems that we will consider in the sequel will be put in a
framework where $D$ is bounded, so that, $\forall t\geq 0$, $D(t)\,\in[0,D_{max}]$; such a bound makes sense in an optimal control
framework since it prevents infinite values of the control, which might occur
when harvesting the photobioreactor. In order to simplify notations, we then
introduce the following change of time and variable
$(t,y)=(D_{max}\tau,\frac{x}{k})$, which yields
\begin{equation}\label{final}
\frac{dy}{dt}=\dot y=\frac{\mu(t)y}{1+y}-ry-uy
\end{equation}
where $r=\frac{\rho}{D_{max}}$ and $u=\frac{D}{D_{max}}\,\in\,[0,1]$ is the
new control. We also have $\mu(t)=\bar \mu=
  \frac{\bar \nu}{\kappa D_{max}}$ for $t\,\in\,[0,\bar T]$ and $0$ for
  $t\,\in\,[\bar T,T ]$ (with $\bar T=D_{max}\bar{\mathcal T}$ and
  $T=D_{max}{\mathcal T}$).

\subsection{Productivity optimization in constant light environment}

In a previous work \cite{Masci2010}, we have studied the productivity
optimization of a microalgae photobioreactor with light-attenuation in the
Droop framework with constant light. In that study, since we wanted to
optimize the long-term productivity, we looked for the control values for $D$
and $s_{in}$ that optimized the instantaneous biomass output flow at
equilibrium, that is
\[
\max_{u} u y^*V
\]
where $V$ is the photobioreactor volume (assumed here to be constant). This study was complex because the shading
was dependent on the internal substrate quota. In the present case, it will greatly
simplify with $s_{in}$ that does not need to be optimized. Indeed, for a given
dilution $u$, the equilibrium of (\ref{final}) in the presence of light is
\[
y^*= \frac{\bar \mu}{r+u}-1
\]
which needs to be non-negative, so that $0\leq u\leq \bar \mu-r$. The
positivity of $u$ imposes that $r\leq \bar \mu$, that is the
respiration needs to be weaker than the maximal growth. For a given $u$, the
productivity rate at equilibrium is then
\[
\frac{\bar \mu u V}{r+u}-uV
\]
whose optimum value is reached in
\begin{equation}\label{usigma}
u_\sigma=\sqrt{\bar\mu r}-r
\end{equation}
which is positive because $r\leq \bar\mu$ but requires
\begin{equation}\label{restr_usigma}
\bar\mu\leq\frac{(r+1)^2}{r}
\end{equation}
to be smaller or equal to $1$  (otherwise, the optimal dilution is $u=1$). This yields the optimal
productivity rate:
\[
(\sqrt{\bar\mu}-\sqrt{ r})^2V
\]
It is important to note that the equilibrium is then
\begin{equation}\label{ysigma}
y_\sigma=\sqrt{\frac{\bar\mu}{r}}-1
\end{equation}
which maximizes the net production rate $\frac{\bar\mu y}{1+y}-ry=uy$. We will use this definition of $y_\sigma$ even when it is not achievable with some $u_\sigma\leq 1$.

\subsection{Productivity optimization in day/night environment}

In an environment with varying light we cannot settle for an instanteneous
productivity rate optimization since this equilibrium cannot be maintained
during the night. In essence, we want to optimize the long term productivity of
the photobioreactor, that is we want that, everyday, the same maximal amount is
produced. The problem that we consider is therefore
\[
\max_{u(t)\,\in\,[0,1]} \int_0^T u(t)V y(t) dt
\]
We then need to add constraints to the solution that we want to obtain;
indeed, at the end of the day, we want to be able to start operating the
photobioreactor in the same conditions for the next day. This then requires that we
add the constraint
\[
y(T)=y(0)
\]

We therefore are faced with the
following optimal control problem
\begin{equation}\label{OCP}
\begin{array}{l}
\displaystyle\max_{u(t)\,\in\,[0,1]} \int_0^T u(t) y(t) dt\\
\hspace*{1cm}\mbox{with}\hspace*{0.5cm}\dot y = \frac{\mu(t) y}{1+y}-ry-uy\\
\hspace*{2.1cm}y(T)=y(0)\\
\end{array}
\end{equation}

\subsubsection{Parameter constraints}

In order to solve this problem, it is convenient to observe that $y(T)=y(0)$
cannot be achieved for large values of $y$ even without considering
optimality. Indeed, for all $t$, we have $\dot y<0$ when  $y(t)>\displaystyle\frac{\bar \mu-r}{r}$ independently of the choice of $u$; therefore, an initial condition such that $y(0)>\displaystyle\frac{\bar \mu-r}{r}$ cannot be considered since necessarily $y(T)<y(0)$ in that case. We then know that, for admissible initial conditions below that threshold, $y(t)$ will stay below this threshold for all times. It also implies that, whenever $u(t)=0$ for such solution with $t\,\in\,[0,\bar T)$, $\dot y>0$ because $y(t)$ then tends toward $\frac{\bar \mu-r}{r}$;

We could make this
bound stronger by noticing that, for a given $y(0)$, the largest value of $y(T)$
that can be achieved is reached by taking $u(t)=0$ for all times; indeed, at any time, applying $u(t)>0$ implies that $\dot y$ is smaller than if $u(t)=0$ were applied. If the value
of $y(T)$ corresponding to $u(t)=0$ is smaller than $y(0)$, then the
corresponding initial condition cannot be part of the optimal
solution. Solving (\ref{final})  with $u(t)=0$ in
the interval $[0,\bar T]$, by separating the variables yields
\[
\frac{r\ln\left(\displaystyle\frac{y(\bar T)}{y_0}\right)-\bar\mu\ln\left(\displaystyle\frac{\bar\mu-r(1+y(\bar
    T))}{\bar\mu-r(1+y_0)}\right)}{r(\bar\mu-r)}=\bar T
\]
where we denoted $y(0)$ as $y_0$. Trivially, the integration of (\ref{final}), for the dark period ($u(t)=0$) on the interval $[\bar T,T]$,  yields
\[
y(T)=y(\bar T)e^{-r(T-\bar T)}
\]
so that, introducing this equation in the previous one, we get
\[
\frac{r\ln\!\left(\!\displaystyle\frac{y(T)e^{r(T-\bar T)}}{y_0}\!\right)\!-\bar\mu\ln\!\left(\!\displaystyle\frac{\bar\mu-r(1+y(T)e^{r(T-\bar T)})}{\bar\mu-r(1+y_0)}\!\right)\!}{r(\bar\mu-r)}\!=\bar T
\]
The equality $y(T)=y_0$ is then achieved with $u(t)=0$ when solving this last
equation for $y_0$ with $y(T)=y_0$, which yields
\[
y_{0max}=\frac{\bar \mu-r}{r}\frac{e^{\frac{r}{\bar \mu}(\bar\mu\bar T-rT)}-1}{e^{\frac{rT}{\bar \mu}(\bar\mu-r)}-1}
\]
For larger values of $y_0$, we have $y(T)<y_0$ independently of the choice of
$u(t)$; for smaller values of $y_0$, there exist control functions $u(t)$ that
guarantee $y(T)=y_0$. The constraint $\mu>r$, which is necessary for growth to
occur in the light phase guarantees that the first fraction and the
denominator of the second one in $y_{0max}$ are positive. We then need to add
the constraint
\begin{equation}\label{mumin}
\bar \mu>\frac{rT}{\bar T}
\end{equation}
to ensure the positivity of $y_{0max}$ and so the possibility of the existence of a solution to the optimal control problem (\ref{OCP}). Note that, in the case where $\bar T=\frac{T}{2}$, this simply
means that $\bar \mu>2r$.

It is also interesting to see that, if a constant control $u(t)=1$ is applied,
a periodic solution is obtained for
\[
y_{0min}=\frac{\bar \mu-r-1}{r+1}\frac{e^{\frac{r+1}{\bar \mu}(\bar\mu\bar T-(r+1)T)}-1}{e^{\frac{(r+1)T}{\bar \mu}(\bar\mu-r-1)}-1}
\]
which can be positive if $\bar \mu>\frac{(r+1)T}{\bar T}$. For any value of
$y_0$ smaller than $y_{0min}$, any control law would force $y(T)>y_0$. As a
consequence, $y_0$, solution of problem (\ref{OCP}), should belong to the
interval $[ y_{0min},y_{0max}]$.

\subsubsection{Maximum principle}

In order to solve problem (\ref{OCP}), we will use Pontryagin's Maximum
Principle (PMP, \cite{Pont}) in looking for a control law maximizing the Hamiltonian
\[
H(x,u,\lambda,t)\triangleq \left[\lambda\left(
\left(\frac{\mu(t)}{1+y}-r\right)y-uy\right)+uy\right]
\]
with the constraint
\[
\left\{\begin{array}{lll}
\dot y&=&\frac{\mu(t)y}{1+y}-ry-uy\\
\dot \lambda&=&\lambda\left(
-\frac{\mu(t)}{(1+y)^2}+r+u\right)-u
\end{array}\right.
\]
In addition, we should add the constraint
\[
\lambda(T)=\lambda(0).
\]
Indeed, the solution of the optimal control problem is
independent of the reference initial time: defining $x(t)=x(t-T)$,  $u(t)=u(t-T)$, and $\lambda(t)=\lambda(t-T)$ for values of $t$ larger than $T$, we have that $x(t)$, $u(t)$  and therefore $\lambda(t)$ are
unchanged if we consider the interval $[t_0,T+t_0]$ (for $0<t_0<T$) rather
than $[0,T]$. Since $\lambda(t)$ is continuous inside the interval when
considering the problem over $[t_0,T+t_0]$, it is continuous in time $T$ and
$\lambda(0)=\lambda(T)$ \cite{Gilbert}.

We see from the form of the Hamiltonian that
\[
\frac{\partial H}{\partial u}=1-\lambda
\]
so that, when $\lambda>1$, we have $u=0$, when $\lambda<1$, we have $u=1$, and
when $\lambda=1$ over some time interval, intermediate singular control is applied.

In the sequel, we propose candidate solutions to the PMP by
making various hypotheses on the value of $\lambda(0)=\lambda_0$.

{\bf Bang-bang with $\lambda_0>1$: } With $\lambda_0>1$, we have $u=0$ at times $0$ and $T$. At any given time $0\leq t \leq \bar T$
before the first switch, the solution of (\ref{final}) yields
\begin{equation}\label{y-phase1}
\frac{r\ln\left(\displaystyle\frac{y(t)}{y_0}\right)-\bar\mu\ln\left(\displaystyle\frac{\bar\mu-r(1+y(t))}{\bar\mu-r(1+y_0)}\right)}{r(\bar\mu-r)}=t
\end{equation}
and, as stated earlier, $y(t)$ is increasing because $y(0)<y_{0max}<\frac{\mu-r}{r}$ .
The constancy of the Hamiltonian during the light phase then imposes that
\begin{equation}\label{lambda-phase1}
\lambda(t) y(t)\left(\frac{\bar\mu}{1+y(t)}-r\right)=\lambda_0 y_0\left(\frac{\bar\mu}{1+y_0}-r\right)
\end{equation}
for all times $t\,\in\,(0,\bar T)$ such that $u(t)=0$. A switch to $1$ then needs to occur
between time $0$ and $T$ (otherwise the payoff would be $0$) and this switch
cannot take place in the dark phase. Indeed, in that zone, as long as
$u(t)=0$, the $\lambda$ dynamics are
\[
\dot \lambda = r\lambda
\]
with $\lambda(t)>1$. The adjoint variable is therefore an increasing function
in that region, and cannot go through $\lambda=1$. We will use this impossibility of switch from $0$ to $1$ in the dark phase several times in the sequel.

For the solution that we
study, a switch then needs to take place at time $t_{01}$ in the $(0,\bar T)$
interval and for $y(t_{01})=y_{01}$ and $\lambda(t_{01})=1$ solutions of (\ref{y-phase1})-(\ref{lambda-phase1}).
\begin{equation}\label{y-switch1}
\frac{r\ln\left(\displaystyle\frac{y_{01}}{y_0}\right)-\bar\mu\ln\left(\displaystyle\frac{\bar\mu-r(1+y_{01})}{\bar\mu-r(1+y_0)}\right)}{r(\bar\mu-r)}=t_{01}
\end{equation}
\begin{equation}\label{lambda-switch1}
y_{01}\left(\frac{\bar\mu}{1+y_{01}}-r\right)=\lambda_0 y_0\left(\frac{\bar\mu}{1+y_0}-r\right)
\end{equation}

Another constraint that appears at the switching instant from $u=0$ to $u=1$
is that $\dot\lambda<0$, which amounts to $\frac{\bar\mu}{(1+y)^2}>r$ or
$y<y_\sigma$ (see (\ref{ysigma})). After time $t_{01}$, $y(t)$ then converges
increasingly or decreasingly toward $\frac{\bar \mu-r-1}{r+1}$

Due to the constancy of the Hamiltonian, another switch can only take place
at time $\tilde t$ before time $\bar T$ if
\[
y(\tilde t)\left(\frac{\bar\mu}{1+y(\tilde t)}-r\right)= y_{01}\left(\frac{\bar\mu}{1+y_{01}}-r\right)
\]
where we have used the fact that $\lambda(\tilde t)=\lambda(t_{01})=1$ at the
switching instants. This can only happen for two values of $y(\tilde t)$:  $y(\tilde t)=y_{01}$ and another value $y(\tilde t) = \frac{\bar\mu}{1+y_{01}}-r$ which is larger than $y_\sigma$. Since
$y(t)$ was converging to $\frac{\bar \mu-r-1}{r+1}$ with $u(t)=1$, $y(\tilde t)$ cannot go through $y_{01}$ again
unless $y_{01}=\frac{\bar \mu-r-1}{r+1}$. In this last case, by considering the $\dot \lambda$ dynamics, we see that another switch could
only take place if $u(t)=1$ solves the conditions for being a singular
solution to the optimal control; this will be handled later. Generically, a
single switch can then only take place inside the interval $(0,\bar
T)$.

The solution then reaches the time $\bar T$ with $(y(t),\lambda(t))=(\bar
y,\bar\lambda)$ that solve the same kind of equations as (\ref{y-phase1}) and (\ref{lambda-phase1}):
\begin{equation}\label{y-Tbar}
\frac{(r+1)\ln\left(\displaystyle\frac{\bar
      y}{y_{01}}\right)-\bar\mu\ln\left(\displaystyle\frac{\bar\mu-(r+1)(1+\bar
      y)}{\bar\mu-(r+1)(1+y_{01})}\right)}{(r+1)(\bar\mu-r-1)}=\bar T-t_{01}
\end{equation}
\begin{equation}\label{lambda-Tbar}
\bar\lambda\bar y\left(\frac{\bar\mu}{1+\bar y}-r-1\right)+\bar y=y_{01}\left(\frac{\bar\mu}{1+y_{01}}-r\right)
\end{equation}
Since $\lambda(\bar T)< 1$ and $\lambda(T)>1$, a switch from $u=1$ to $u=0$
then needs to take place inside the $(\bar T,T)$ interval. With the dynamics being in the form
\[
\dot y=-(r+1)y\hspace*{1cm}\dot \lambda =(r+1)\lambda-1
\]
another switch can only take place if $\bar\lambda>\frac{1}{r+1}$; otherwise
$\lambda$ cannot go through $1$ again. The switching point $(t_{10},y_{10})$
is then characterized by
\begin{equation}\label{y-switch2}
y_{10}= \displaystyle\bar ye^{-(r+1)(t_{10}-\bar T)}
\end{equation}
\begin{equation}\label{lambda-switch2}
\lambda(t_{10})=1= \displaystyle\bar\lambda e^{(r+1)(t_{10}-\bar T)}-\displaystyle\frac{e^{(r+1)(t_{10}-\bar T)}-1}{r+1}
\end{equation}
After this switching, the dynamics become
\[
\dot y=-ry\hspace*{1cm}\dot \lambda =r\lambda
\]
so that no other switch can take place and these dynamics and the constraints
$y(T)=y_0$ and $\lambda(T)=\lambda_0$ impose that
\begin{equation}\label{y-switch3}
y_0= \displaystyle y_{10}e^{-r(T-t_{10})}
\end{equation}
\begin{equation}\label{lambda-switch3}
\lambda_0= \displaystyle e^{r(T-t_{10})}
\end{equation}
In the end, we have a system of 8 algebraic equations\linebreak
(\ref{y-switch1})-(\ref{lambda-switch3}) with eight unknowns, which we
solve numerically.

Even though, we were not able to lead this study analytically all the way to
the end, we have shown the qualitative form of the solutions analytically. It
is made of four phases:
\begin{itemize}
\item Growth with a closed photobioreactor until a sufficient biomass level is
  reached
\item Maximal harvesting of the photobioreactor with simultaneous growth
\item Maximal harvesting of the photobioreactor with no growth until a low level of
  biomass is reached
\item Passive photobioreactor: no harvesting, no growth, only respiration
\end{itemize}
The first two phases take place in the presence of light, the other two in the
dark. In phase 3, harvesting of as much biomass  produced in the light
phase as possible is continued while not going below the level where  the residual biomass left is
sufficient to efficiently start again the next day.

{\bf Bang-singular-bang with $\lambda_0>1$: }

We will first look at what a singular arc could be. For that, we see that
$\frac{\partial H}{\partial u}=1-\lambda$ should be 0 over a time interval and
compute its time derivatives.
\[
\frac{d}{dt}\left(\frac{\partial H}{\partial u}\right)_{\arrowvert \lambda=1}=-\frac{\mu(t)}{(1+y)^2}+r
\]
When $\mu(t)=0$, that is in the dark phase, no singular arc is thus
possible. When $\mu(t)=\bar \mu$, this derivative is equal to zero when
$y=y_\sigma$ defined in (\ref{ysigma}). The singular control is then the
control that maintains this equilibrium, that is $u_\sigma=\sqrt{\bar\mu r}-r$
defined in (\ref{usigma}). This control is positive thanks to (\ref{mumin})
but it is smaller or equal to $1$ only if
\begin{equation}\label{mumax}
\bar\mu\leq \frac{(r+1)^2}{r}
\end{equation}
No singular control can exist otherwise. When a singular branch appears in the
optimal solution, it is locally optimal because the second order Kelley condition
$$
\frac{\partial}{\partial u}\left(\frac{d^2}{d\tau^2}\frac{\partial H}{\partial
    u}\right)=
\frac{2\lambda\mu}{(1+y)^3}\geq 0
$$
is satisfied on the singular arc \cite{Kelley}.

The construction of the solution is
very similar to that in the purely bang-bang case. Similarly, a switch needs
to  occur in the interval $(0,\bar T)$. This switch can be from $u=0$ to $u=1$
or from $u=0$ to $u=u_\sigma$ and should occur with $y\leq y_\sigma$ in order
to have $\dot\lambda\leq 0$. In fact, if a switch first occurs to $u=1$, an
argument identical to the one in the previous section shows that no switch
back to $0$ can take place before $\bar T$; this same argument can in fact be
used to show that no switch to $u=u_{\sigma}$ can take place either since: in
both cases, $\lambda$ should get back to $1$, which we show to be impossible.

A switch from $0$ to $u_\sigma$ then takes place once $\lambda=1$ at
$(t_{0\sigma},y_\sigma)$. Equations (\ref{y-phase1})-(\ref{lambda-phase1}) can
then be used to identify this switching instant:
\begin{equation}\label{y-switchsigma}
\frac{r\ln\left(\displaystyle\frac{y_{\sigma}}{y_0}\right)-\bar\mu\ln\left(\displaystyle\frac{\bar\mu-r(1+y_{\sigma})}{\bar\mu-r(1+y_0)}\right)}{r(\bar\mu-r)}=t_{0\sigma}
\end{equation}
\begin{equation}\label{lambda-switchsigma}
y_{\sigma}\left(\frac{\bar\mu}{1+y_{\sigma}}-r\right)=\lambda_0 y_0\left(\frac{\bar\mu}{1+y_0}-r\right)
\end{equation}
From there, $\lambda(t)=1$ and $y(t)=y_{\sigma}$ for some time. This could be
until $t=\bar T$, followed directly by $u=0$ in the dark phase but, more
generically, the singular arc ends at time $t_{\sigma 1}<\bar T$, where a switch
occurs toward $u=1$. From then on, things
are unchanged with respect to the bang-bang case. The equations that define
the transitions from $t_{\sigma 1}$ to $\bar T$ are similar to (\ref{y-Tbar})
and (\ref{lambda-Tbar}):
 \begin{equation}\label{y-Tbar-sigma}
\frac{(r+1)\ln\left(\displaystyle\frac{\bar
      y}{y_\sigma}\right)-\bar\mu\ln\left(\displaystyle\frac{\bar\mu-((r+1)(1+\bar
      y)}{\bar\mu-(r+1)(1+y_\sigma)}\right)}{(r+1)(\bar\mu-r-1)}=\bar T-t_{\sigma 1}
\end{equation}
\begin{equation}\label{lambda-Tbar-sigma}
\bar\lambda\bar y\left(\frac{\bar\mu}{1+\bar y}-r-1\right)+\bar y=y_\sigma\left(\frac{\bar\mu}{1+y_\sigma}-r\right)
\end{equation}
The remainder of the solution is unchanged with respect to the bang-bang one,
so that we can compute the solution by solving system
(\ref{y-switch2})-(\ref{lambda-switch3}) and
(\ref{y-switchsigma})-(\ref{lambda-Tbar-sigma}) of eight algebraic equations
with eight unknown variables.

Again, the analytical approach has helped us identify the qualitative form of
the optimal productivity solution. It now contains five phases:
\begin{itemize}
\item Growth with a closed photobioreactor until a sufficient biomass level is
  reached
\item Maximal equilibrium productivity rate on the singular arc
\item Maximal harvesting of the photobioreactor with simultaneous growth
\item Maximal harvesting of the photobioreactor with no growth until a low level of
  biomass is reached
\item Passive photobioreactor: no harvesting, no growth, only respiration
\end{itemize}

For this form of solution, we see that maximal instantaneous productivity is
achieved during the whole second phase, when the singular solution occurs.

{\bf Solution with $\lambda_0<1$: }

Such a solution would mean that harvesting takes place during the whole dark
phase because no transition from $u=0$ to $u=1$ can take place in this
phase, as we have already shown. Two possibilities then occur: either $u=1$ all the time or switches
from $u=1$ to $u=0$ or $u_\sigma$ and then back to $u=1$ take place in the
interval $(0,\bar T)$.

In the latter case, the first switch from  $u=1$ to $u=0$ can only take place
with $y>y_\sigma$ because of the constraint that $\dot\lambda>0$ with
$\lambda=1$ at that moment. Then, when the control $u=0$ is applied
for some time, the solution $y(t)$ is increasing. We also have that the switch from $u=0$
to $u=1$ can only take place
with $y<y_\sigma$ because of the constraint that $\dot\lambda<0$ with
$\lambda=1$ at that moment. This is in contradiction with the fact that $y(t)$
was increasing from above $y_\sigma$.

We can also show that no strategy in the $(0,\bar T)$ interval can have the
form $u=1\rightarrow u_\sigma \rightarrow u=0 \text{ or } 1$. Indeed, in order to reach the
singular arc with $u=1$, a solution should be coming from above it. If the
switch that takes place at the end of the singular phase is from $u_\sigma$ to
$0$, $y(t)$ will increase and there should be a subsequent switch from $0$ to
$1$ which is impossible with $y(t)>y_\sigma$. If the switch that takes place
at the end of the singular phase is from $u_\sigma$ to $1$, $y(t)$ will
decrease all the time between $t_{\sigma 1}$ and $T$, which is in
contradiction with the fact that we had $y(0)>y_\sigma$.

\begin{figure}[t]
\centerline{\includegraphics[width=.85\columnwidth]{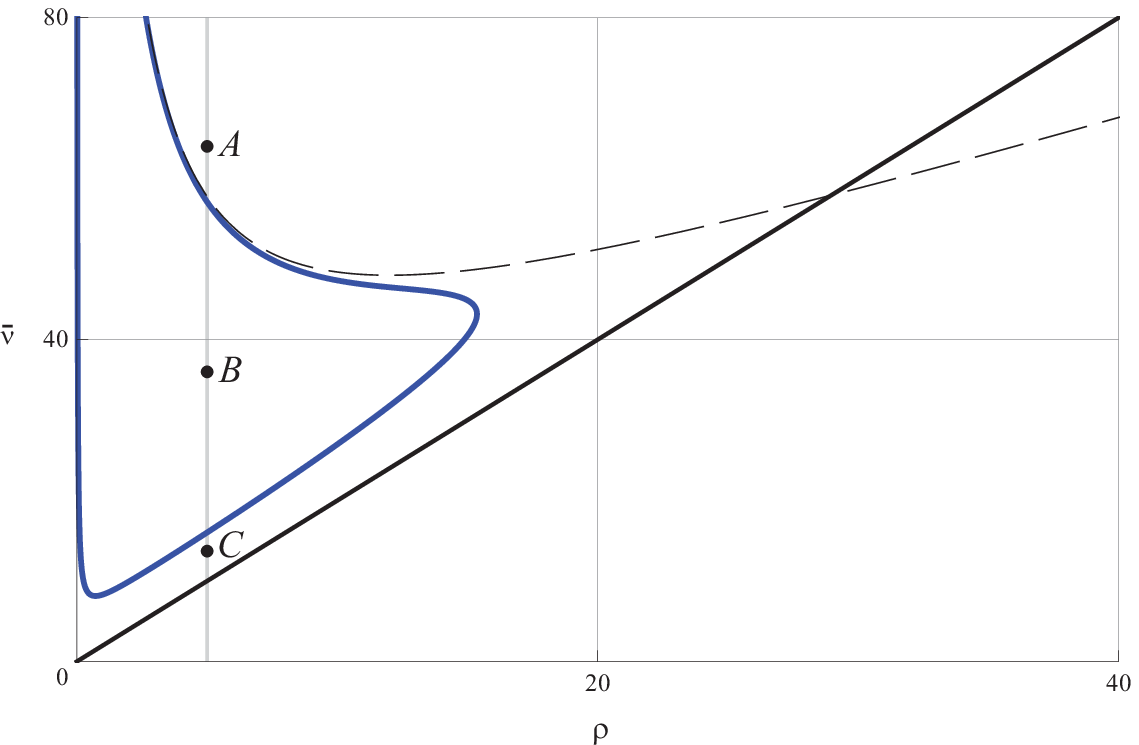}}
\caption{Bifurcation picture for $D_{max}=12$, $\kappa=1$, $T=1$, $\bar T=T/2$. The solid black line is $\bar\nu={\kappa\rho T}/{\bar T}$ (see (\ref{mumin})), the dashed line is $\bar\nu = \kappa(\rho+D_{max})^2/\rho$ and it is related to (\ref{restr_usigma}). Optimal patterns for $A$, $B$ and $C$ are shown on Fig.~3}\label{fig3}
\bigskip
\centerline{\includegraphics[width=.75\columnwidth]{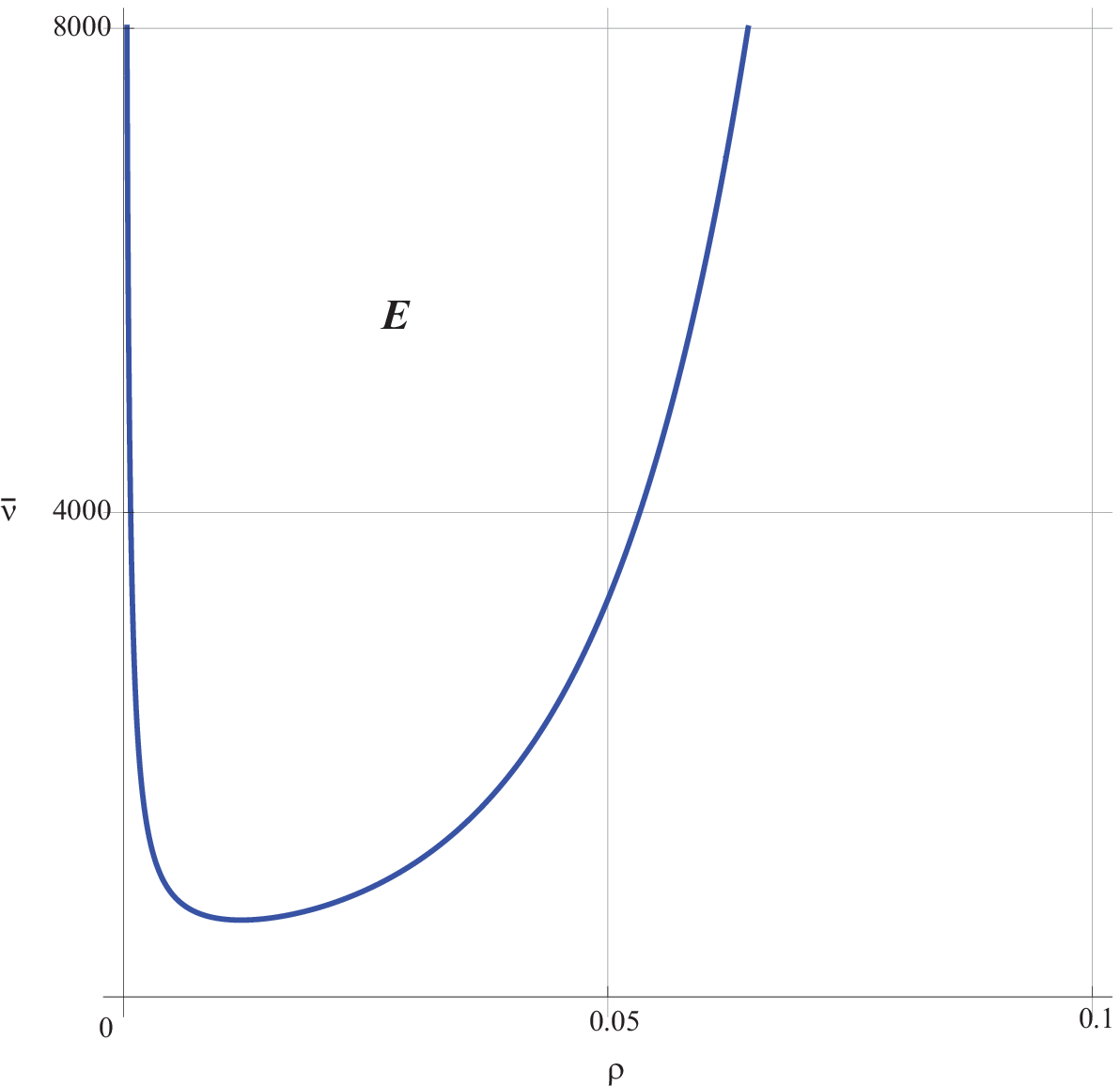}}
\caption{The region $E$ with optimal constant control $u=1$. Below this region, this diagram is connected with Fig.~1}\label{fig3}
\end{figure}

The only potential optimal control in that family is therefore $u(t)=1$ for
all times. Using the expressions computed previously, this control can be a
candidate optimal control law only if $y_0=y_{0min}$ as we have seen earlier
and the complete dynamics should satisfy:
\begin{equation}\label{lambda-Tbar-1}
\bar\lambda\bar y\left(\frac{\bar\mu}{1+\bar y}-r-1\right)+\bar y=\lambda_0 y_{0}\left(\frac{\bar\mu}{1+y_{0}}-r\right)
\end{equation}
\begin{equation}\label{y-1}
\bar y= \displaystyle y_{0}e^{(r+1)(T-\bar T)}
\end{equation}
\begin{equation}\label{lambda-1}
\lambda_0= \displaystyle\bar\lambda e^{(r+1)(T-\bar T)}-\displaystyle\frac{e^{(r+1)(T-\bar T)}-1}{r+1}
\end{equation}
with $\lambda_0<1$ and $\bar\lambda<1$.

\section{Bifurcation analysis}

\begin{figure}[t]
\centerline{\includegraphics[width=.8\columnwidth]{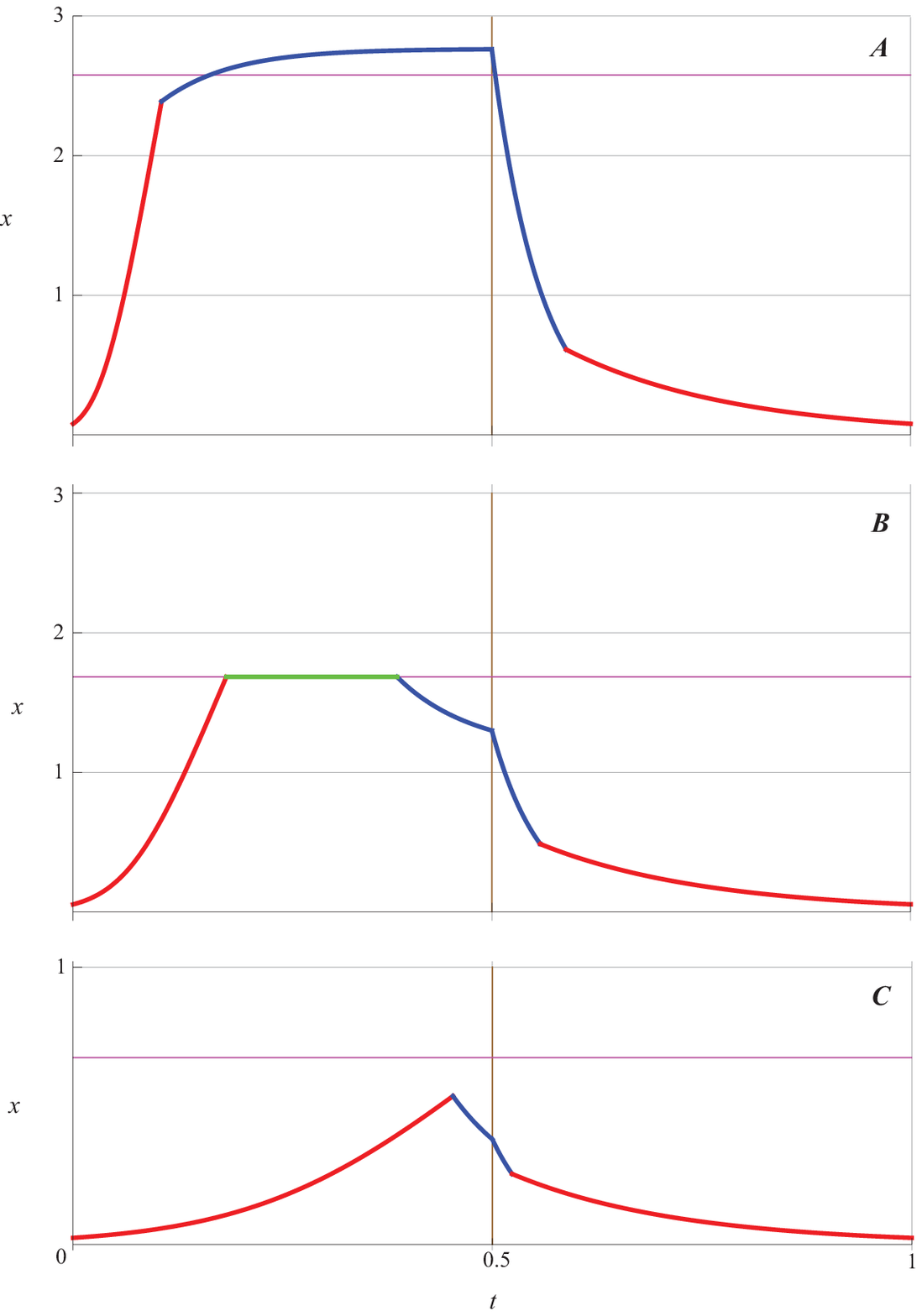}}
\caption{Different optimal patterns: (A) $\bar\nu=14$, (B) $\bar\nu=36$, (C) $\bar\nu=64$; $\rho=5$, $D_{max}=12$, $\kappa=1$, $T=1$, $\bar T=T/2$. Red: $u=0$, Blue:  $u=1$, Green: intermediate control $u\in[0,1]$. Magenta line indicates the level $x=\kappa y_\sigma$ (see (\ref{ysigma}))}\label{fig1}
\end{figure}

In this section, we will consider fixed values of all parameters except of $\bar\nu$ and $\rho$. We build a bifurcation diagram for these two parameters by identifying in which region no solution is possible (where (\ref{mumin}) is not satisfied, it is below the solid black line on Fig.~1), and where the optimal solution is  bang-singular-bang (Fig.~1, inside the blue curve), bang-bang (Fig.~1, outside the blue curve and above the solid black line), and constant at value $1$ (see Fig.~2). But the last case is only realized for extremely large values of $\bar\nu$. We see that the region where singular control can exist is smaller than what is defined by condition (\ref{mumax}). This is due to the fact that, though the singular control is possible, there is not enough time for the control to reach that level (see Fig.~3(C)). For larger values of $\bar\nu$, no singular control is possible and the optimal solution in the light region goes toward the equilibrium corresponding to $u=1$ (see Fig.~3(A)). In that case, as well as in the bang-singular-bang case, the solutions go to the optimal solution of the constant light problem (Fig.~3(B)).

\section{Conclusions}

We have shown that, because of the day-night constraint, the productivity rate cannot be as high as it could have been without it. However, when the maximal growth rate is sufficiently larger than the respiration rate, we manage to have a temporary phase where the productivity rate is at or near this level. The maximal harvesting at the end of the light phase and at the beginning of the dark phase minimizes the biomass during the dark phase and, consequently, the net respiration.
If the maximal growth rate is very large, the optimal solution consists in constantly applying maximal control because the biomass that is built-up in the light phase needs to be harvested even during the night.

\end{document}